[1994/12/01]
\documentclass{ijmart-mod}
\chardef\bslash=`\\ 

\usepackage{graphicx}
\usepackage[breaklinks=true]{hyperref}
\usepackage[ruled,linesnumbered]{algorithm2e}

\newtheorem{thm}{Theorem}[section]
\newtheorem{cor}[thm]{Corollary}
\newtheorem{lem}[thm]{Lemma}
\newtheorem{prop}[thm]{Proposition}

\theoremstyle{definition}
\newtheorem{defn}[thm]{Definition}
\newtheorem{rem}[thm]{Remark}

\theoremstyle{remark}

\newcommand{\eval}[2][\right]{\relax
  \ifx#1\right\relax \left.\fi#2#1\rvert}

\begin{document}
\title[An oracle for zonotope computation]{A linear optimization oracle for~zonotope~computation}

\author[Antoine Deza]{Antoine Deza}
\address{McMaster University, Hamilton, Ontario, Canada}
\email{deza@mcmaster.ca} 

\author[Lionel Pournin]{Lionel Pournin}
\address{LIPN, Universit{\'e} Paris 13, Villetaneuse, France}
\email{lionel.pournin@univ-paris13.fr} 

\begin{abstract}
A class of counting problems ask for the number of regions of a central hyperplane arrangement. By duality, this is the same as counting the vertices of a zonotope. We give several efficient algorithms, based on a linear optimization oracle, that solve this and related enumeration problems. More precisely, our algorithms compute the vertices of a zonotope from the set of its generators and inversely, recover the generators of a zonotope from its vertices. A variation of the latter algorithm also allows to decide whether a polytope, given as its vertex set, is a zonotope and when it is not, to compute its greatest zonotopal summand.
\end{abstract}
\maketitle

\section{Introduction}\label{CZ.sec.introduction}

Linear optimization consists in finding a vertex of a polyhedron that maximizes some linear functional. It is widely used in many areas of science and engineering. Linear optimization is instrumental in the solution to some prominent questions and has led to formulating a number of other widely studied problems.
While linear optimization itself is known to be polynomial time solvable in the size of the problem (see for instance \cite{IllesTerlaky2002} and references therein), the complexity of simplex methods, pivot-based linear optimization algorithms, is still not known. One of the questions that arises from the study of linear optimization, in relation with this particular class of algorithms, is that of the largest diameter a polyhedron can have. Here, by the diameter of a polyhedron, we mean the diameter of the graph made up of its vertices and edges. The largest possible diameter of a polyhedron has been studied from a number of different perspectives \cite{BorgwardtDeLoeraFinhold,DeLoeraHemmeckeLee2015,DeLoeraKaferSanita2019,KitaharaMizuno2013,Sanita2018,TerlakyZhang1993,Ye2011}, and in particular as a function of its dimension and number of facets \cite{KalaiKleitman1992,KleeWalkup1967,Naddef1989,Santos2012,Sukegawa2019}, two parameters that reflect the number of variables and the number of constraints of a linear optimization problem. In practice, the vertices of polyhedra often have rational coordinates and, up to the multiplication by an integer, these vertices are contained in the integer lattice. This is a reason why the diameter of lattice polytopes, the polytopes whose coordinates of vertices are integers, is also widely studied \cite{DelPiaMichini2016,DezaManoussakisOnn2018,DezaPournin2018,KleinschmidtOnn1992,Naddef1989}. In this case, the largest possible diameter is estimated in terms of the dimension and the size of the smallest hypercube the polytope is contained in. \phantom{\cite{Ziegler1995}}

In contrast to linear optimization that consists in finding just one vertex of a polyhedron, convex hull computation amounts to enumerate all the faces of a polytope. A number of efficient algorithms have been given that address this particular problem \cite{AvisBremnerSeidel1997,AvisFukuda1992,PreparataShamos1985}. Since the number of faces of a polytope of arbitrary dimension is exponential (for instance in its number of vertices), the worst case complexity of these algorithms is exponential.

This article treats a case that lies in between linear optimization and convex hull computation. Just as with linear optimization, we are interested in the small dimensional faces of polytopes (their vertices, but also their edges) and, as convex hull computations do, one of our main goals is to enumerate them. We are going to do that for a particular class of polytopes, the zonotopes or, in other words, the Minkowski sums of line segments. By their definition, zonotopes are much more combinatorial, and sometimes behave very differently than arbitrary polytopes. For instance, linear optimization on a zonotope is linear time solvable in the number of its generators. Zonotopes arise in a number of counting problems related to very different fields, often in terms of their dual hyperplane arrangement \cite{BlockWeinberger1992,DezaManoussakisOnn2018,GutekunstMeszarosPetersen2019,MelamedOnn2014}. These counting problems ask about the number of vertices of a zonotope, where the zonotope itself is given as the set of its generators. 
The first contribution of this article is an efficient, convex hull free algorithm that solves this kind of problems in practice. In other words, it enumerates the vertices of a zonotope from the set of its generators. The complexity of this algorithm is linear in the number of vertices of the zonotope and polynomial in the number of generators. Our second contribution is an efficient algorithm that performs the inverse computation, also without carrying out any convex hull. Given the vertex set of a polytope $P$, this algorithm will decide if $P$ is a zonotope, and in this case it will return its set of generators. In this sense it can also be considered a decision algorithm. It is polynomial in the number of vertices of the considered polytope. 

\begin{figure}
\begin{centering}
\includegraphics{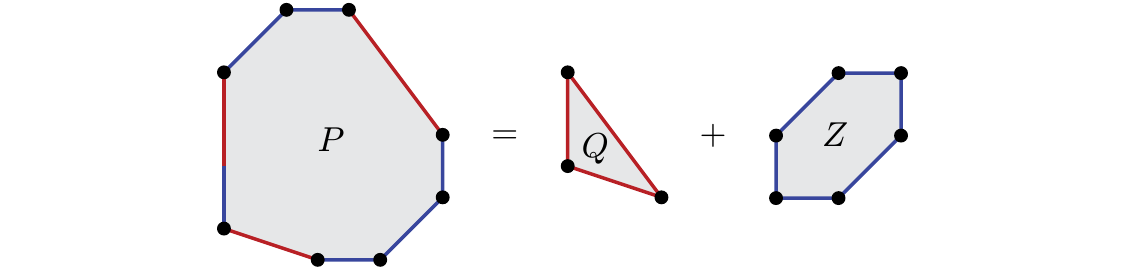}
\caption{The octagon $P$ is the Minkowski sum of $Q$ and $Z$.}\label{Fig.CZ.1}
\end{centering}
\end{figure}

We also provide a third algorithm, that provides a practical take on the question of polytope decomposability, another topic that has attracted significant attention \cite{DezaPournin2019,Kallay1982,Meyer1974,PrzeslawskiYost2008,Shephard1963}. This algorithm is an intermediate step towards our algorithm that computes the generators of a zonotope from its vertex set. It efficiently computes the \emph{greatest zonotopal summand} of an arbitrary polytope. Let us illustrate this notion. By a summand of a polytope $P$, we mean any polytope $Q$ such that $P$ is the Minkowski sum of $Q$ with another polytope. 
In Figure \ref{Fig.CZ.1} for instance, the octagon $P$ is the Minkowski sum of the triangle $Q$ with the hexagon $Z$. In particular, $Q$ and $Z$ are two summands of $P$. 
Observe that $Z$ is, up to translation, the Minkowski sum of three of its edges. In other words, $Z$ is a zonotope.
On the other hand, since $Q$ is a triangle, no line segment---and therefore no zonotope---can be a summand of $Q$. In this case, $Z$ is what we call the greatest zonotopal summand of $P$.

When a polytope $P$ is given as the convex hull fo a finite set of points, linear optimization is polynomial time solvable in the number of these points since it amounts to compute the value of a linear map at each of them. According to the theory developed in \cite{GrotschelLovaszSchrijver1993}, deciding whether a given point is a vertex of $P$ is then also polynomial time solvable. Similarly, deciding whether two points are the extremities of an edge of $P$ can be done in polynomial time. Our algorithms rely on the ability to solve these two problems. In Section \ref{CZ.sec.oracle}, we will give an explicit linear optimization oracle that provides a practical way to do that. In section \ref{CZ.sec.properties}, we will recall a number of properties of Minkowski sums and zonotopes, and derive other properties that will be used in the sequel. The algorithm that enumerates the vertices of a zonotope from its set of generators is given in Section~\ref{CZ.sec.computingZ}. The greatest zonotopal summand of a polytope is defined and studied at in Section~\ref{CZ.sec.computingGZS}, and the algorithm that computes it is described at the end of the section. Finally, the algorithm that enumerates the generators of a zonotope from its vertex set is given in Section \ref{CZ.sec.decidingZ}.

%
%
%
%
%
%
%

\section{A linear optimization oracle}\label{CZ.sec.oracle}

We begin the section with a linear optimization oracle that allows to tell whether the convex hull of a finite subset of $\mathbb{R}^d$ is disjoint from the affine hull of another finite subset of $\mathbb{R}^d$. We then show that this oracle provides a practical way to decide when a polytope is a face of another when both are given as convex hulls of finite sets of points. As a consequence, we obtain an explicit algorithm that efficiently computes the graph of a polytope given either as the set of its vertices or as the convex hull of a finite subset of $\mathbb{R}^d$. At the end of the section, we show how our oracle also allows to compute the rays of a pointed cone given as the conic hull of a set of points.

Consider a finite subset $\mathcal{A}$ of $\mathbb{R}^d$ and a subset $\mathcal{F}$ of $\mathcal{A}$. The convex hull of $\mathcal{A}\mathord{\setminus}\mathcal{F}$ and the affine hull of $\mathcal{F}$ are non-disjoint if and only a convex combination of $\mathcal{A}\mathord{\setminus}\mathcal{F}$ coincides with an affine combination of $\mathcal{F}$; that is, if and only if there exists a family $(\alpha_a)_{a\in\mathcal{A}}$ of real numbers such that

\begin{alignat}{4}
\label{CZ.oracle.eq.1}
\sum_{a\in\mathcal{A}\mathord{\setminus}\mathcal{F}}\alpha_aa & \;- & \,\sum_{a\in\mathcal{F}}\alpha_aa & = 0\mbox{,}\\
\label{CZ.oracle.eq.2}
\sum_{a\in\mathcal{A}\mathord{\setminus}\mathcal{F}}\alpha_a\;\; & & & =1\mbox{,}\\
\label{CZ.oracle.eq.3}
 & & \sum_{a\in\mathcal{F}}\alpha_a & =1\mbox{,}
\end{alignat}
and $\alpha_a\geq0$ when $a\in\mathcal{A}\mathord{\setminus}\mathcal{F}$. In other words, checking whether $\mathrm{conv}(\mathcal{A}\mathord{\setminus}\mathcal{F})$ and $\mathrm{aff}(\mathcal{F})$ are non-disjoint amounts to find a solution to a system of $d+2$ linear equalities (note here that (\ref{CZ.oracle.eq.1}) accounts for $d$ of these equalities) and $|\mathcal{A}|-|\mathcal{F}|$ linear inequalities that state the non-negativity of some of the variables. This feasibility problem, which we denote by $(LO_{\mathcal{A},\mathcal{F}})$ is polynomial time solvable in $|\mathcal{A}|$, $d$, and the binary size $L$ of the input (see~\cite{GrotschelLovaszSchrijver1993} or \cite{IllesTerlaky2002}). In our case and throughout the section, $L$ is the number of bits needed to store $\mathcal{A}$.

Let us explain how $(LO_{\mathcal{A},\mathcal{F}})$ allows for an efficient way to compute the graph of a polytope given as the convex hull of a finite subset of $\mathbb{R}^d$.

\begin{prop}\label{CZ.oracle.prop.1}
Consider a finite subset $\mathcal{A}$ of $\mathbb{R}^d$. A point $x$ in $\mathcal{A}$ is a vertex of $\mathrm{conv}(\mathcal{A})$ if and only if it does not belong to $\mathrm{conv}(\mathcal{A}\mathord{\setminus}\{x\})$.
\end{prop}
\begin{proof}
The vertices of a polytope are precisely its extreme points. Therefore, a point $x$ in $\mathcal{A}$ is not a vertex of $\mathrm{conv}(\mathcal{A})$ if and only if it can be written as a convex combination of $\mathcal{A}$ whose coefficient for $x$ is less than $1$. Note that $x$ appears on both sides of that equality. Solving this equality for $x$ results in an equivalent equation that expresses $x$ as a convex combination of $\mathcal{A}\mathord{\setminus}\{x\}$.
\end{proof}

Recall that the affine hull of a single point only contains that point. Hence, according to Proposition \ref{CZ.oracle.prop.1}, A point $x$ in $\mathcal{A}$ is a vertex of $\mathrm{conv}(\mathcal{A})$ if and only if $(LO_{\mathcal{A},\{x\}})$ is not feasible. Solving this feasability problem for every point $x$ in $\mathcal{A}$ allows to recover the vertex set of $\mathrm{conv}(\mathcal{A})$ in polynomial time in $|\mathcal{A}|$, $d$, and $L$. Proposition \ref{CZ.oracle.prop.1} can be generalized as follows.

\begin{prop}\label{CZ.oracle.prop.2}
Consider a finite subset $\mathcal{A}$ of $\mathbb{R}^d$. The convex hull of a subset $\mathcal{F}$ of $\mathcal{A}$ is a face of the convex hull of $\mathcal{A}$ if and only if the affine hull of $\mathcal{F}$ is disjoint from the convex hull of $\mathcal{A}\mathord{\setminus}\mathcal{F}$.
\end{prop}
\begin{proof}
Let $\mathcal{F}$ be a subset of $\mathcal{A}$. Consider the orthogonal projection
$$
\pi:\mathbb{R}^d\rightarrow\mathrm{aff}(\mathcal{F})^\perp\mbox{,}
$$
where $\mathrm{aff}(\mathcal{F})^\perp$ the orthogonal complement of $\mathrm{aff}(\mathcal{F})$ in $\mathbb{R}^d$; that is, the set of the vectors in $\mathbb{R}^d$ orthogonal to it. By construction, $\pi$ sends all the points in $\mathcal{F}$ to a single point $x$. Now observe that $\mathrm{aff}(\mathcal{F})$ is disjoint from $\mathrm{conv}(\mathcal{A}\mathord{\setminus}\mathcal{F})$ if and only if $x$ is not contained in the convex hull of $\pi(\mathcal{A}\mathord{\setminus}\mathcal{F})$. By Proposition \ref{CZ.oracle.prop.1}, this is equivalent to $x$ being a vertex of the convex hull of $\pi(\mathcal{A})$ which, in turn, is equivalent to $\mathrm{conv}(\mathcal{F})$ being a face of $\mathrm{conv}(\mathcal{A})$.
\end{proof}

Consider two distinct elements $x$ and $y$ of $\mathcal{A}$. According to Proposition \ref{CZ.oracle.prop.2}, $\mathrm{conv}(\{x,y\})$ is an edge of $\mathrm{conv}(\mathcal{A})$ if and only if $(LO_{\mathcal{A},\{x,y\}})$ is not feasible. Therefore, in order to compute the graph of the convex hull of a finite subset $\mathcal{A}$ of $\mathbb{R}^d$, it is sufficient to solve $(LO_{\mathcal{A},\{x\}})$ for every point $x$ in $\mathcal{A}$ in order to recover the vertex set $\mathcal{V}$ of $\mathrm{conv}(\mathcal{A})$ and then to solve $(LO_{\mathcal{A},\{x,y\}})$ for any pair of distinct points $x$ and $y$ in $\mathcal{V}$ in order to decide whether $\mathrm{conv}(\{x,y\})$ is an edge of $\mathrm{conv}(\mathcal{A})$. Since the number of feasibility problems to solve is quadratic in $|\mathcal{A}|$ and each of these problems is polynomial time solvable in $|\mathcal{A}|$, $d$, and $L$, we immediately obtain the following complexity result.

\begin{thm}\label{CZ.oracle.thm.1}
Consider a finite subset $\mathcal{A}$ of $\mathbb{R}^d$. The graph of $\mathrm{conv}(\mathcal{A})$ can be computed in polynomial time in $|\mathcal{A}|$, $d$, and $L$.
\end{thm}

Note that the whole $k$-skeleton of $\mathrm{conv}(\mathcal{A})$ cannot be computed in polynomial time using the same ideas when $k$ is greater than $1$. In this case, the number of feasibility problems to solve is exponential in $|\mathcal{A}|$. Indeed, the possible candidates for being the vertex set of a $2$-dimensional face of $\mathrm{conv}(\mathcal{A})$ would be all the sets of at least three vertices of $\mathrm{conv}(\mathcal{A})$. However, for any given positive integer $k$, it is possible to compute in polynomial time all the faces of $\mathrm{conv}(\mathcal{A})$ with at most $k$ vertices and arbitrary dimension. According to Proposition \ref{CZ.oracle.prop.2}, this amounts to solve $(LO_{\mathcal{A},\mathcal{F}})$ for every non-empty subset $\mathcal{F}$ of $\mathcal{A}$ with at most $k$ elements, whose number is at most a degree $k$ polynomial in $|\mathcal{A}|$. 

We now explain how the same ideas allow to compute in polynomial time the rays of a pointed cone spanned by a finite subset $\mathcal{A}$ of $\mathbb{R}^d\mathord{\setminus}\{0\}$. This will be used in Section~\ref{CZ.sec.computingZ} in order to compute the vertices of a zonotope efficiently from its generators. Recall that the \emph{cone spanned by $\mathcal{A}$}, or \emph{conic hull of $\mathcal{A}$}, is the polyhedral cone made up of all the linear combinations with non-negative coefficients of the points in $\mathcal{A}$. This cone, which we denote by $\mathrm{cone}(\mathcal{A})$, is pointed when it admits $0$ as a vertex. Note that, when $\mathcal{A}$ is made up of a single non-zero point, $\mathrm{cone}(\mathcal{A})$ is a half-line incident to $0$.
%
%
%

Note that, if $\mathcal{A}$ contains a pair of linearly dependent points, then the cone spanned by $\mathcal{A}$ is not affected if one of these points is removed from $\mathcal{A}$. Hence, we can assume that any two points in $\mathcal{A}$ are linearly independent. The following proposition is the conic counterpart to Proposition \ref{CZ.oracle.prop.1}.

\begin{prop}\label{CZ.oracle.prop.3}
Consider a finite subset $\mathcal{A}$ of pairwise linearly independent points of $\mathbb{R}^d\mathord{\setminus}\{0\}$ that spans a pointed cone. The half-line spanned by a point $x$ in $\mathcal{A}$ is a ray of the cone spanned by $\mathcal{A}$ if and only if the line through $0$ and $x$ is disjoint from the convex hull of $\mathcal{A}\mathord{\setminus}\{x\}$.
\end{prop}
\begin{proof}
Consider a point $x$ in $\mathcal{A}$ and assume that the half-line $\mathrm{cone}\{x\}$ is a ray of $\mathrm{cone}(\mathcal{A})$. Let $H$ be a supporting hyperplane of $\mathrm{cone}(\mathcal{A})$ such that
$$
\mathrm{cone}(\mathcal{A})\cap{H}=\mathrm{cone}\{x\}\mbox{.}
$$

Since any two points in $\mathcal{A}$ are linearly independent, $x$ is the only element of $\mathcal{A}$ contained in $\mathrm{cone}\{x\}$ and, therefore, in $H$. As a consequence, $\mathcal{A}\mathord{\setminus}\{x\}$ is contained in one of the open half-spaces bounded by $H$ and its convex hull is necessarily disjoint from the line through $0$ and $x$.

Now assume that the line through $0$ and $x$ is disjoint from the convex hull of $\mathcal{A}\mathord{\setminus}\{x\}$. By Proposition \ref{CZ.oracle.prop.2}, $\mathrm{conv}\{0,x\}$ is an edge of $\mathrm{conv}(\mathcal{A}\cup\{0\})$. Consider a supporting hyperplane $H$ of $\mathrm{conv}(\mathcal{A}\cup\{0\})$ such that
$$
\mathrm{conv}(\mathcal{A}\cup\{0\})\cap{H}=\mathrm{conv}\{0,x\}\mbox{.}
$$

Since $\mathrm{conv}(\mathcal{A}\cup\{0\})\mathord{\setminus}\mathrm{conv}\{0,x\}$ is contained in one of the open half-spaces bounded by $H$, and since $H$ contains $0$, then $\mathrm{cone}(\mathcal{A})\mathord{\setminus}\mathrm{cone}(\{x\})$ is also contained in that half-space. Hence, $H$ is a supporting hyperplane of the cone spanned by $\mathcal{A}$ and it intersects this cone along the half-line spanned by $x$. In other words, that half-line is a ray of the cone spanned by $\mathcal{A}$.
\end{proof}

Now observe that the line through $0$ and a point $x$ in $\mathcal{A}$ is disjoint from the convex hull of $\mathcal{A}\mathord{\setminus}\{x\}$ if and only if $(LO_{\mathcal{A}\cup\{0\},\{0,x\}})$ is not feasible. Again, this feasibility problem is polynomial time solvable in $|\mathcal{A}|$, $d$, and $L$. In particular, it follows from Proposition \ref{CZ.oracle.prop.3} that solving this problem allows to tell whether the half-line spanned by $x$ is a ray of the conic hull of $\mathcal{A}$.

By these observations, we obtain the following.

\begin{thm}\label{CZ.oracle.thm.2}
Consider a finite subset $\mathcal{A}$ of pairwise linearly independent points of $\mathbb{R}^d\mathord{\setminus}\{0\}$ that spans a pointed cone. The rays of $\mathrm{cone}(\mathcal{A})$ can be computed in polynomial time in $|\mathcal{A}|$, $d$, and $L$.
\end{thm}

The input of some of the algorithms we describe in the sequel are polytopes given as the set $\mathcal{V}$ of their vertices. In fact, these polytopes could also be given as any finite set $\mathcal{A}$ of points they are the convex hull of. In this case, the complexity of these algorithms would be in terms of $|\mathcal{A}|$ instead of $|\mathcal{V}|$.

\section{Combinatorial properties of Minkowski additions}\label{CZ.sec.properties}


Recall that the Minkowski sum of two subsets $P$ and $Q$ of $\mathbb{R}^d$ is
$$
P+Q=\{x+y:(x,y)\in{P\mathord{\times}Q}\}\mbox{.}
$$

When $P$ and $Q$ are polyhedra, the faces of $P+Q$ are exactly the Minkowski sums of a face $F$ of $P$ and a face $G$ of $Q$ such that, for some non-zero vector $c$ in $\mathbb{R}^d$, the linear functional $x\mapsto{c\mathord{\cdot}x}$ is maximized at $F$ in $P$ and at $G$ in $Q$ (see for instance Proposition 12.1 in~\cite{Fukuda2015} or Lemma~2.1 in \cite{DezaPournin2019}). As already mentioned, a zonotope, is the Minkowski sum of a finite set of line segments. In fact, a zonotope $Z$ contained in $\mathbb{R}^d$ is uniquely obtained, up to translation, as the Minkowski sum of a finite set of pairwise non-homothetic line segments incident to $0$ and whose first non-zero coordinate of the other vertex is positive. We refer to these particular line segments as the \emph{generators} of $Z$.

Now recall that a summand of a polytope $P$ is a polytope $Q$ such that $P$ is obtained as the Minkowski sum of $Q$ with another polytope.
\begin{figure}[b]
\begin{centering}
\includegraphics{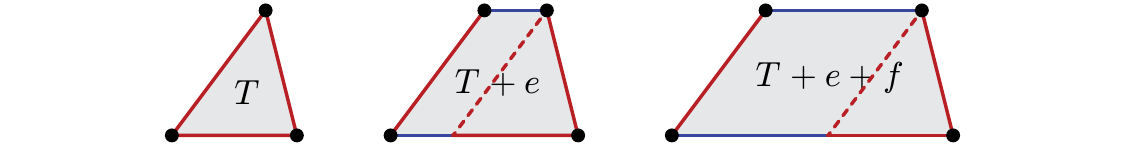}
\caption{A triangle $T$, its Minkowski sum with a horizontal line segment $e$ (center), and the Minkowski sum of $T+e$ with a line segment $f$ homothetic to $e$ (right).}\label{Fig.CZ.2}
\end{centering}
\end{figure}
We borrow the following decomposability characterization from \cite{DezaPournin2019} that we will use to compute the $1$-dimensional summands of a polytope efficiently. 

\begin{lem}[{\cite[Theorem 2.5]{DezaPournin2019}}]\label{CZ.properties.lem.1}
A polytope $P$ has a summand homothetic to a polytope $Q$ if and only if $P$ and $P+Q$ have the same number of vertices.
\end{lem}

Lemma \ref{CZ.properties.lem.1} is illustrated in Figure \ref{Fig.CZ.2} in the case when $P$ is the quadrilateral $T+e$, that admits a line segment $e$ as a summand, and $Q$ is a line segment $f$ homothetic to $e$. As can be seen, $T+e+f$ is still a quadrilateral. 

We also introduce counting arguments that will be used to speedup the computations. If $P$ is a polytope and $s$ is a line segment, we denote by $\langle{s}\rangle_P$ the set of the edges of $P$ homothetic to $s$ or, equivalently, parallel to $s$.

\begin{lem}\label{CZ.properties.lem.2}
Consider a  $d$-dimensional polytope $P$ and a line segment $s$, both contained in $\mathbb{R}^d$. If $s$ is a summand of $P$, then
\begin{itemize}
\item[(i)] At least $d$ edges of $P$ are contained in $\langle{s}\rangle_P$,
\item[(ii)] All the edges of $P$ in $\langle{s}\rangle_P$ are at least as long as $s$,
\item[(iii)] Every shortest element of $\langle{s}\rangle_P$ is a summand of $P$.
\end{itemize}
\end{lem}
\begin{proof}
Assume that there exists a polytope $Q$ such that $P$ is the Minkowski sum of $Q$ and $s$ and consider the orthogonal projection
$$
\pi:\mathbb{R}^d\rightarrow\mathrm{aff}(s)^\perp\mbox{.}
$$

Note that the image by $\pi$ of any element of $\langle{s}\rangle_P$ is a vertex of $\pi(P)$. In fact, $\pi$ induces a bijection between $\langle{s}\rangle_P$ and the vertex set of $\pi(Z)$. Indeed, first observe that, by convexity, $\pi$ cannot send distinct elements of $\langle{s}\rangle_P$ to the same point. Further consider a vertex $v$ of $\pi(P)$ and a vector $c$ contained in $\mathrm{aff}(s)^\perp$ such that the map $x\mapsto{c\mathord{\cdot}x}$ is maximized at $v$ in $\pi(P)$. Observe that $\pi(P)$ coincides, up to translation, with $\pi(Q)$. Hence, $x\mapsto{c\mathord{\cdot}x}$ is maximized at a vertex of $\pi(Q)$ and, therefore, at a vertex or an edge of $Q$ parallel to $s$. Denote by $Q_v$ this vertex or edge of $Q$.  Since $c$ is contained in $\mathrm{aff}(s)^\perp$, the map $x\mapsto{c\mathord{\cdot}x}$ is constant within $s$. Hence, according to Lemma~2.1 in \cite{DezaPournin2019},
\begin{equation}\label{CZ.properties.lem.2.eq.1}
\pi^{-1}(\{v\})\cap{P}=Q_v+s\mbox{.}
\end{equation}

In other words, the set of the points in $P$ that $\pi$ sends to $v$ is an edge of $P$ obtained as the Minkowski sum of $Q_v$ and $s$. Therefore, this edge belongs to $\langle{s}\rangle_P$ and cannot be shorter than $s$. In particular assertion (ii) holds.

Since $\pi(Q)$ is a polytope of dimension $d-1$, it has at least $d$ vertices. As $\pi$ projects distinct segments in $\langle{s}\rangle_P$ to distinct points, $\langle{s}\rangle_P$ cannot contain less than $d$ elements. In other words, assertion (i) holds.

Now observe that, if $Q_v$ is an edge of $Q$ for every vertex $v$ of $\pi(P)$, then $Q$ still admits a summand homothetic to $s$. Hence, $P$ has a summand homothetic to $s$ that is longer than $s$. We can assume without loss of generality that $s$ is a longest such summand of $P$. In this case, $Q_v$ is a vertex of $Q$ for some vertex $v$ of $\pi(P)$. By Equality (\ref{CZ.properties.lem.2.eq.1}), some edge of $P$ is a translate of $s$. According to assertion (ii), all the shortest edges of $P$ in $\langle{s}\rangle_P$ must be translates of $s$ and, therefore, summands of $P$, which proves assertion (iii).
\end{proof}

In the case of zonotopes, the statement of Lemma \ref{CZ.properties.lem.2} can be refined. In order to do that, we will use the following property of zonotopes.

\begin{prop}\label{CZ.properties.prop.1}
A $d$-dimensional zonotope has at least $2^d$ vertices.
\end{prop}
\begin{proof}
Consider a zonotope $Z$ contained in $\mathbb{R}^d$. If $Z$ is $d$-dimensional, then it admits $d$ generators that are not all contained in a hyperplane of $\mathbb{R}^d$. The Minkowski sum of these generators is a $d$-dimensional combinatorial hypercube $C$. By construction, $C$ is a summand of $Z$. Therefore, according to Lemma 2.3 from \cite{DezaPournin2019}, there exists an injection from the vertex set of $C$ into the vertex set of $Z$. Since $C$ has $2^d$ vertices, the proposition follows.
\end{proof}

Lemma \ref{CZ.properties.lem.2} can be improved as follows in the case of zonotope.

\begin{lem}\label{CZ.properties.lem.3}
If $Z$ is a $d$-dimensional zonotope then, for every generator $g$ of $Z$, $\left|\langle{g}\rangle_Z\right|\!\geq{2^{d-1}}$ and all the elements of $\langle{g}\rangle_Z$ are translates of $g$.
\end{lem}
\begin{proof}
Consider a $d$-dimensional polytope $Z$ contained in $\mathbb{R}^d$ and a generator $g$ of $Z$. By the definition, two distinct generators of $Z$ cannot be parallel. Since the faces of a Minkowski sum of polytopes is a Minkowski sum of their faces, every edge of $Z$ is a translate of one of its generators. As a consequence, every element of $\langle{g}\rangle_Z$ is necessarily a translate of $g$.

Consider the orthogonal projection
$$
\pi:\mathbb{R}^d\rightarrow\mathrm{aff}(g)^\perp\mbox{.}
$$

Observe that $\pi(Z)$ is a zonotope of dimension $d-1$ obtained, up to translation, as the Minkowski sum of the images by $\pi$ of the generators of $\pi(Z)$. 
By the argument in the proof of Lemma~\ref{CZ.properties.lem.2}, $\pi$ induces a bijection between $\langle{g}\rangle_Z$ and the vertex set of $\pi(Z)$. Since $\pi(Z)$ is a $(d-1)$-dimensional zonotope, the desired result therefore follows from Proposition \ref{CZ.properties.prop.1}.
\end{proof}

By Lemma \ref{CZ.properties.lem.3}, $\left|\langle{S}\rangle_Z\right|$ is at least $2^{d-1}$ when $S$ is a generator of a $d$-dimensional zonotope Z. Note that $\left|\langle{S}\rangle_Z\right|$ is not necessarily a multiple of $2^{d-1}$ in general. For instance, the rhombic dodecahedron is a $3$-dimensional zonotope whose exactly $6$ edges are obtained as translates of each generator.

\section{An efficient algorithm to compute zonotopes}\label{CZ.sec.computingZ}

Throughout this section, $Z$ is a fixed zonotope contained in $\mathbb{R}^d$. Recall that a zonotope is, up to translation, the Minkowski sum of its generators. In this section, we assume that $Z$ is exactly the Minkowski sum of its generators, which can be done without loss of generality by translating $Z$, if needed. Denote by $\mathcal{G}$ the set of the non-zero vertices of the generators of $Z$. The purpose of the section is to give an algorithm to recover the vertex set of $Z$ from $\mathcal{G}$. 
Observe that, since $Z$ is the Minkowski sum of its generators, it is also equal to the convex hull of all the possible subsums of $\mathcal{G}$; that is,
$$
Z=\mathrm{conv}\left\{\sum_{x\in\mathcal{X}}x:\mathcal{X}\subset\mathcal{G}\right\}\mbox{,}
$$
where, by convention, the sum of the elements of $\mathcal{X}$ is equal to $0$ when $\mathcal{X}$ is empty. In particular, every vertex of $Z$ is the sum of a unique subset of $\mathcal{G}$. However, not all the subsets of $\mathcal{G}$ sum to a vertex of $Z$. The following theorem characterizes the subsets of $\mathcal{G}$ that have this property.


\begin{thm}\label{CZ.computingZ.thm.1}
Consider a subset $\mathcal{X}$ of $\mathcal{G}$. The sum of the points in $\mathcal{X}$ is a vertex of $Z$ if and only if the set $[-\mathcal{X}]\cup[\mathcal{G}\mathord{\setminus}\mathcal{X}]$ spans a pointed cone.
\end{thm}
\begin{proof}
Consider the zonotope $Z_\mathcal{X}$ equal to the Minkowski sum of the line segments incident to $0$ and whose other vertex is a point in $[-\mathcal{X}]\cup[\mathcal{G}\mathord{\setminus}\mathcal{X}]$. Note that $Z$ and $Z_\mathcal{X}$ are translates of one another. More precisely,
$$
Z=Z_\mathcal{X}+x\mbox{,}
$$
where $x$ denotes the sum of the points in $\mathcal{X}$. As a consequence, $x$ is a vertex of $Z$ if and only if $0$ is a vertex of $Z_\mathcal{X}$. Since $Z_\mathcal{X}$ is the Minkowski sum of its generators, it is equal to the convex hull of all the possible subsums of $[-\mathcal{X}]\cup[\mathcal{G}\mathord{\setminus}\mathcal{X}]$. In particular $Z_\mathcal{X}$ admits $[-\mathcal{X}]\cup[\mathcal{G}\mathord{\setminus}\mathcal{X}]$ as a subset. Hence, if $0$ is a vertex of $Z_{\mathcal{X}}$, then $[-\mathcal{X}]\cup[\mathcal{G}\mathord{\setminus}\mathcal{X}]$ is contained in one of the open half-spaces limited by a hyperplane through $0$ and therefore spans a pointed cone.

Since $Z_\mathcal{X}$ is the convex hull of all the subsums of $[-\mathcal{X}]\cup[\mathcal{G}\mathord{\setminus}\mathcal{X}]$, it must contain $0$ and be contained in the cone spanned by $[-\mathcal{X}]\cup[\mathcal{G}\mathord{\setminus}\mathcal{X}]$. Hence, if $[-\mathcal{X}]\cup[\mathcal{G}\mathord{\setminus}\mathcal{X}]$ spans a pointed cone, then $0$ is a vertex of $Z_\mathcal{X}$.
\end{proof}

Let us illustrate Theorem \ref{CZ.computingZ.thm.1} by showing that $0$ is a vertex of $Z$. Since the first non-zero coordinate of every point in $\mathcal{G}$ is positive, $0$ is not a convex combination of $\mathcal{G}$. In this case, according to Proposition \ref{CZ.oracle.prop.1}, the convex hull of $\mathcal{G}\cup\{0\}$ admits $0$ as a vertex. As a consequence, the cone spanned by $\mathcal{G}$ is pointed and, by Theorem \ref{CZ.computingZ.thm.1}, $0$ is a vertex of $Z$.

It is worth noting that the condition in Theorem \ref{CZ.computingZ.thm.1} can be checked efficiently. More precisely, the cone spanned by $[-\mathcal{X}]\cup[\mathcal{G}\mathord{\setminus}\mathcal{X}]$ is pointed if and only if the following system of $|\mathcal{G}|$ linear inequalities is feasible.
\begin{alignat}{2}
\label{CZ.computingZ.eq.1}
c\mathord{\cdot}g\geq1 & \mbox{ for all }g\in\mathcal{X}\mbox{,}\\
\label{CZ.computingZ.eq.2}
-c\mathord{\cdot}g\geq1 & \mbox{ for all }g\in\mathcal{G}\mathord{\setminus}\mathcal{X}\mbox{.}
\end{alignat}

Indeed, the feasibility of this system is equivalent to the existence of a vector $c\in\mathbb{R}^d$ such that the map $x\mapsto{c\mathord{\cdot}x}$ is maximized exactly at $0$ within the cone spanned by $[-\mathcal{X}]\cup[\mathcal{G}\mathord{\setminus}\mathcal{X}]$. In other words, this cone admits $0$ as a vertex.

Theorem \ref{CZ.computingZ.thm.1} already provides a way to compute the vertices of~$Z$. Indeed, in order to do that, it suffices to check, for each subset $\mathcal{X}$ of $\mathcal{G}$ whether $0$ is contained in the convex hull of $[-\mathcal{X}]\cup[\mathcal{G}\mathord{\setminus}\mathcal{X}]$. It is possible though, that many of these subsets do not sum to a vertex of $Z$. In order to avoid considering these subsets, we will use the following lemma.

\begin{lem}\label{CZ.computingZ.lem.1}
Consider a subset $\mathcal{X}$ of $\mathcal{G}$. If $\mathcal{X}$ sums to a vertex $x$ of $Z$, then the vertices of $Z$ adjacent to $x$ are the sums of $x$ with any element of $[-\mathcal{X}]\cup[\mathcal{G}\mathord{\setminus}\mathcal{X}]$ that spans a ray of the conic hull of $[-\mathcal{X}]\cup[\mathcal{G}\mathord{\setminus}\mathcal{X}]$.
\end{lem}
\begin{proof}
Assume that the sum of the points in $\mathcal{X}$ is a vertex $x$ of $Z$. In this case, according to Theorem \ref{CZ.computingZ.thm.1}, $[-\mathcal{X}]\cup[\mathcal{G}\mathord{\setminus}\mathcal{X}]$ spans a pointed cone. As in the proof of Theorem \ref{CZ.computingZ.thm.1}, we consider the zonotope $Z_\mathcal{X}$ whose generators are incident to $0$ on one end and to a point in $[-\mathcal{X}]\cup[\mathcal{G}\mathord{\setminus}\mathcal{X}]$ on the other. This zonotope is a translate of $Z$. More precisely,
$$
Z=Z_\mathcal{X}+x\mbox{.}
$$

According to this, in order to prove the lemma, we only need to show that the vertices of $Z_\mathcal{X}$ adjacent to $0$ are exactly the points in $[-\mathcal{X}]\cup[\mathcal{G}\mathord{\setminus}\mathcal{X}]$ that span a ray of the conic hull of $[-\mathcal{X}]\cup[\mathcal{G}\mathord{\setminus}\mathcal{X}]$.

By construction, $Z_\mathcal{X}$ is the Minkowski sum of its generators. Hence, $Z_\mathcal{X}$ is the convex hull of all the possible subsums of $[-\mathcal{X}]\cup[\mathcal{G}\mathord{\setminus}\mathcal{X}]$. In particular, it is contained in the cone spanned by $[-\mathcal{X}]\cup[\mathcal{G}\mathord{\setminus}\mathcal{X}]$. Now recall that the edges of a zonotope are translates of its generators. Therefore, the vertices of $Z_\mathcal{X}$ adjacent to $0$ must be among the points from $[-\mathcal{X}]\cup[\mathcal{G}\mathord{\setminus}\mathcal{X}]$. Consider a point $y$ in $[-\mathcal{X}]\cup[\mathcal{G}\mathord{\setminus}\mathcal{X}]$. The segment with vertices $0$ and $y$ is an edge of $Z_\mathcal{X}$ if and only if there exists a supporting hyperplane $H$ of $Z_\mathcal{X}$ such that
$$
Z_\mathcal{X}\cap{H}=\mathrm{conv}\{0,y\}\mbox{.}
$$

Since the points of $\mathrm{cone}([-\mathcal{X}]\cup[\mathcal{G}\mathord{\setminus}\mathcal{X}])$ are precisely the multiples by a non-negative coefficient of the points in $Z_\mathcal{X}$, this is equivalent to
$$
\mathrm{cone}([-\mathcal{X}]\cup[\mathcal{G}\mathord{\setminus}\mathcal{X}])\cap{H}=\mathrm{cone}\{y\}\mbox{.}
$$

In other words, $\mathrm{cone}\{y\}$ is a ray of the cone spanned by $[-\mathcal{X}]\cup[\mathcal{G}\mathord{\setminus}\mathcal{X}]$.
\end{proof}

According to Theorem \ref{CZ.oracle.thm.2}, the condition in the statement of Lemma \ref{CZ.computingZ.lem.1} can be checked efficiently using the oracle described in Section \ref{CZ.sec.oracle}. As discussed above, this condition can also be checked by solving the feasibility problem made up of the inequalities (\ref{CZ.computingZ.eq.1}) and (\ref{CZ.computingZ.eq.2}) for each generator $g$ of $Z$, where $\mathcal{X}$ is replaced by $\mathcal{X}\cup\{g\}$ if $g$ does not belong to $\mathcal{X}$ and by $\mathcal{G}\mathord{\setminus}\{g\}$ otherwise.

Let us now give an informal description of our algorithm that computes the vertices of $Z$ from $\mathcal{G}$. Recall that $Z$ admits $0$ as a vertex. Our algorithm starts from that vertex and computes all the vertices of $Z$ adjacent to it. According to Lemmas \ref{CZ.oracle.thm.2} and \ref{CZ.computingZ.lem.1}, this can be done in polynomial time in $|\mathcal{G}|$, $d$, and the binary size $L$ required to store $\mathcal{G}$. Then the procedure is repeated greedily and computes the vertices of $Z$ adjacent to the new vertices of $Z$ that have been discovered, and so on until the neighbors of all the discovered vertices have been computed. Since the graph of a polytope (made up of its vertices and edges) is connected, this indeed computes all the vertices of $Z$. In order to further speedup our algorithm, we use the following proposition that allows to compute only a subset of the neighbors of each vertex.

\begin{prop}\label{CZ.computingZ.prop.1}
If $\mathcal{X}$ is a non-empty subset of $\mathcal{G}$ that sums to a vertex $x$ of $Z$, then there exists a vertex $y$ of $Z$ adjacent to $x$ such that $x-y\in\mathcal{X}$.
\end{prop}
\begin{proof}
The proof proceeds by induction on the dimension of $Z$. If $Z$ has dimension $0$, then $\mathcal{G}$ is empty and the desired result immediately holds. Assume that the dimension of $Z$ is positive, and that the desired statement holds for any zonotope of dimension less than the dimension of $Z$.

Consider a non-empty subset $\mathcal{X}$ of $\mathcal{G}$ that sums to a vertex $x$ of $Z$. Since $\mathcal{X}$ is non-empty and the first non-zero coordinate of any point it contains is positive, then $x$ is necessarily distinct from $0$ and its first non-zero coordinate must be positive. We will review two cases.

Assume that the first coordinate of $x$ is positive. As $Z$ contains $0$, $Z$ must have an edge incident to $x$ whose first coordinate of the other vertex is less than that of $x$. By Lemma \ref{CZ.computingZ.lem.1}, the other vertex $y$ of this edge is such that either $x-y$ belongs to $\mathcal{X}$ or to $-[\mathcal{G}\mathord{\setminus}\mathcal{X}]$. Since the first coordinate of the points in $\mathcal{G}$ is non-negative, $x-y$ must belong to $\mathcal{X}$, as desired.

Now assume that the first coordinate of $x$ is equal to $0$. Since the first coordinate of the points in $\mathcal{G}$ is non-negative, it follows that the first coordinate of all the points in $\mathcal{X}$ must be equal to $0$. Moreover, the zonotope generated by the elements of $\mathcal{G}$ whose first coordinate is equal to $0$ is a proper face of $Z$. Hence, the proposition holds by induction.
\end{proof}

\begin{algorithm}[b]\label{CZ.computingZ.algorithm.1}

$\mathcal{T}\leftarrow\{0\}$

$\xi(0)\leftarrow\emptyset$

$\mathcal{V}\leftarrow\emptyset$

\While{$\mathcal{T}\neq\emptyset$}{
Pick an element $x$ of $\mathcal{T}$

\For{every point $y$ in $x+\mathcal{G}\mathord{\setminus}\xi(x)$}{
\If{$y$ does not belong to $\mathcal{T}\cup\mathcal{V}$}{
\If{$y-x$ spans a ray of $\mathrm{cone}([-\xi(x)]\cup[\mathcal{G}\mathord{\setminus}\xi(x)])$}{
$\mathcal{T}\leftarrow\mathcal{T}\cup\{y\}$ and $\xi(y)\leftarrow\xi(x)\cup\{y-x\}$
}
}
}
$\mathcal{T}\leftarrow\mathcal{T}\mathord{\setminus}\{x\}$ and $\mathcal{V}\leftarrow\mathcal{V}\cup\{x\}$
}

Return $\mathcal{V}$
\caption{Computing $\mathcal{V}$ from $\mathcal{G}$}
\end{algorithm}

Consider a vertex $x$ of $Z$ and the subset $\mathcal{X}$ of $\mathcal{G}$ it is the sum of. A consequence of Proposition \ref{CZ.computingZ.prop.1} is that $x$ can be reached from $0$ by a path in the graph of $Z$ that visits only vertices equal to subsums of $\mathcal{X}$. In other words, in order to discover new vertices of $Z$ from $x$ in the algorithm sketched above, one only needs to check the points $y$ such that $y-x$ belongs to $\mathcal{G}\mathord{\setminus}\mathcal{X}$, and the algorithm will still compute all the vertices of $Z$. This is what Algorithm \ref{CZ.computingZ.algorithm.1} does.

Let us give a detailed description of Algorithm \ref{CZ.computingZ.algorithm.1}. In this algorithm, $\mathcal{T}$ is the set of the vertices of $Z$ that have been discovered, but not treated yet in the sense that their neighbors in the graph of $Z$ are still to be computed. The set of the vertices that have been treated, in the same sense is denoted by $\mathcal{V}$. Initially, $\mathcal{T}$ only contains $0$ and $\mathcal{V}$ is empty. Upon completion of the algorithm, $\mathcal{V}$ is the set of the vertices of $Z$. For each point $x$ in $\mathcal{T}\cup\mathcal{V}$, the subset of $\mathcal{G}$ that sums to $x$ is denoted by $\xi(x)$. For instance, $\xi(0)$ is equal to $\emptyset$.

While $\mathcal{T}$ is non-empty, the algorithm picks a point $x$ from $\mathcal{T}$, and considers all the points $y$ that are the sum of $x$ with an element of $\mathcal{G}\mathord{\setminus}\xi(x)$. By Proposition~\ref{CZ.computingZ.prop.1}, one can restrict to only consider these points to enumerate the vertex set of $Z$. In Line $7$, the algorithm first checks whether $y$ has not been discovered yet (which can be done in logarithmic time in the number of vertices of $Z$ using an appropriate data structure). If $y$ has not been discovered, the algorithm checks in Line~$9$ whether $y$ is a vertex of $Z$, using the condition stated by Lemma~\ref{CZ.computingZ.lem.1}. According to Theorem \ref{CZ.oracle.thm.2}, this can be done in polynomial time in $|\mathcal{G}|$, $d$, and the binary size $L$ required to store $\mathcal{G}$. If $y$ is a vertex of $Z$, then it is inserted in $\mathcal{T}$ and $\xi(y)$ is computed in Line $10$. Once $x$ has been treated, it is removed from $\mathcal{T}$ and placed in $\mathcal{V}$ in Line $14$.
%

\begin{thm}
There exists a polynomial function $p:\mathbb{R}^3\rightarrow\mathbb{R}$ such that the vertex set of a $d$-dimensional zonotope with $n$ vertices and $m$ generators can be computed from the set of its generators in time $O(n\,p(m,d,L))$, where $L$ is the number of bits required to store all these generators.
\end{thm}
\begin{proof}
By Theorem \ref{CZ.oracle.thm.2}, there exists a polynomial function $q:\mathbb{R}^3\rightarrow\mathbb{R}$ such that the test in Line~$9$ of Algorithm \ref{CZ.computingZ.algorithm.1} can be done in time $O(q(m,d,L))$ for a $d$-dimensional zonotope with $m$ generators, where $L$ is the binary size required to store these generators. Hence, according to the description of the algorithm, the vertex set of a $d$-dimensional zonotope with $n$ vertices and $m$ generators can be computed from the set of its generators in time
$$
O(nm[q(m,d,L)+\mathrm{log}\,n])\mbox{.}
$$

Since $n$ is at most $2^m$, we obtain the desired result.
\end{proof}

\section{The greatest zonotopal summand of a polytope}\label{CZ.sec.computingGZS}

We introduce the greatest zonotopal summand of a polytope in this section. We also discuss some of its properties and give an efficient algorithm to compute it for a polytope given as the set of its vertices. In the remainder of the section $P$ is a fixed $d$-dimensional polytope with $n$ vertices. 

Denote by $\mathcal{E}$ the set made up of the edges of $P$ that are also summands of~$P$ and consider a segment $e$ in $\mathcal{E}$. 
We refer to as $e^\star$ the unique translate of $e$ whose one vertex is incident to $0$ and whose first non-zero coordinate of the other vertex is positive. According to Lemma~\ref{CZ.properties.lem.2}, any edge $f$ of $P$ in the intersection $\mathcal{E}\cap\langle{E}\rangle_P$ has the same length as $e$ and, therefore, $e^\star$ and $f^\star$ must coincide. In the remainder of the section, we consider the set
$$
\mathcal{G}=\{e^\star:e\in\mathcal{E}\}\mbox{.}
$$

While in Section \ref{CZ.sec.computingZ}, $\mathcal{G}$ was a set of points, here $\mathcal{G}$ contains line segments. However, in both cases, $\mathcal{G}$ describes the generators of a zonotope.

\begin{defn}\label{CZ.computingGZS.defn.1}
We call \emph{greatest zonotopal summand of $P$} and denote by $z(P)$ the Minkowski sum of the line segments contained in $\mathcal{G}$.
\end{defn}

By this definition, $z(P)$ is a zonotope. Let us show that $z(P)$ is indeed a summand of $P$ and that it is the greatest such summand.

\begin{thm}\label{CZ.computingGZS.thm.1}
There exists a polytope $r(P)$ with no $1$-dimensional summand such that the Minkowski sum $z(P)+r(P)$ is equal to $P$.
\end{thm}
\begin{proof}
Note that distinct line segments in $\mathcal{G}$ cannot be parallel. Hence, none of these segments admit a summand homothetic to another. Since each line segment in $\mathcal{G}$ is a summand of $P$, their Minkowski sum is necessarily a summand of $P$. Therefore, there exists a polytope $r(P)$ such that
\begin{equation}\label{CZ.computingGZS.thm.1.eq.1}
P=z(P)+r(P)\mbox{.}
\end{equation}

Now assume that $r(P)$ has a $1$-dimensional summand $s$. In this case, $s$ is also a summand of $P$ and, according to Lemma \ref{CZ.properties.lem.2}, so are the shortest elements of $\langle{s}\rangle_P$. Let $e$ be a shortest element of $\langle{s}\rangle_P$. By construction, $e^\star$ must belong to $\mathcal{G}$ and is therefore a generator of $z(P)$. According to (\ref{CZ.computingGZS.thm.1.eq.1}), $e^\star+s$ is then a summand of $P$. However, by Lemma \ref{CZ.properties.lem.2}, $e^\star+s$ should be shorter than $e$, proving that $s$ cannot be a summand of $r(P)$ in the first place.
\end{proof}

\begin{cor}\label{CZ.computingGZS.cor.1}
If $P$ admits a zonotope $Z$ as a summand, then $z(P)$ necessarily also admits $Z$ as a summand.
\end{cor}
\begin{proof}
By Theorem \ref{CZ.computingGZS.thm.1},
$$
P=z(P)+r(P)\mbox{,}
$$
where $r(P)$ does not admit a $1$-dimensional summand. Therefore, if a zonotope is a summand of $P$, then all of its generators must be summands of $z(P)$. Hence, that zonotope is itself a summand of $z(P)$.
\end{proof}

Recall that $n$ denotes the number of vertices of $P$. The edges of $P$ are quadratically-many in $n$ and, as shown in Section \ref{CZ.sec.oracle} they can all be computed in polynomial time in $n$, $d$, and the binary size $L$ required to store all the vertices of $P$. It turns out that the vertex set of the Minkowski sum of $P$ with a line segment can also be computed in polynomial time in $n$, $d$, and $L$. In fact, we have the following more general observation.

\begin{rem}
Consider two finite subsets $\mathcal{A}$ and $\mathcal{B}$ of $\mathbb{R}^d$. According to Theorem~\ref{CZ.oracle.thm.1}, the graph of the Minkowski sum $\mathrm{conv}(\mathcal{A})+\mathrm{conv}(\mathcal{B})$ can be computed in polynomial time in $|\mathcal{A}||\mathcal{B}|$, $d$, and the number of bits required to store the points in $\mathcal{A}$ and in $\mathcal{B}$. Indeed, this amounts to compute the graph the convex hull of the $\mathcal{A}+\mathcal{B}$, a subset of at most $|\mathcal{A}||\mathcal{B}|$ points of $\mathbb{R}^d$. 
\end{rem}
\begin{algorithm}[t]\label{CZ.computingGZS.algorithm.1}
$\mathcal{R}\leftarrow\emptyset$

\For{every subset $\{x,y\}$ of $\mathcal{V}$ such that $x\neq{y}$}{
$E\leftarrow\mathrm{conv}\{x,y\}$

\If{$e$ is an edge of $P$}{
\eIf{some segment $s$ in $\mathcal{R}$ is parallel to $e$}{
$\mu(s)\leftarrow\mu(s)+1$

\If{$e$ is shorter than $s$}{
$s\leftarrow{e^\star}$
}
}{
$\mathcal{R}\leftarrow\mathcal{R}\cup\{e^\star\}$ and $\mu(e^\star)\leftarrow1$
}
}
}

$\mathcal{W}\leftarrow\mathcal{V}$ and $\mathcal{G}\leftarrow\emptyset$

\For{every segment $s$ in $\mathcal{R}$}{
\If{$\mu(s)\geq{d}$}{
Compute the vertex set $\mathcal{X}$ of $\mathrm{conv}(\mathcal{W})+s$

\If{$|\mathcal{X}|=|\mathcal{W}|$}{
$\mathcal{G}\leftarrow\mathcal{G}\cup\{s\}$

\For{every point $x$ in $\mathcal{W}\mathord{\setminus}\mathcal{X}$}{
Replace $x$ in $\mathcal{W}$ by $x-s^+$
}
}
}
}

Return $\mathcal{G}$ and $\mathcal{W}$
\caption{Computing $\mathcal{G}$ and $\mathcal{W}$ from $\mathcal{V}$}
\end{algorithm}

As a consequence of this remark, the set $\mathcal{G}$ of the generators of $\mathrm{z}(P)$ can be computed in polynomial time in $n$, $d$, and $L$ as well. Algorithm \ref{CZ.computingGZS.algorithm.1} is a polynomial time algorithm in $n$, $d$, and $L$ that not only computes $\mathcal{G}$, but also the vertex set of $r(P)$. This vertex set is denoted by $\mathcal{W}$ in the algorithm. The vertex set of $P$, denoted by $\mathcal{V}$, is the only input of the algorithm. Algorithm \ref{CZ.computingGZS.algorithm.1} is split in two parts. The first part, from Line $1$ to Line $14$ computes a set $\mathcal{R}$ of candidates for belonging to $\mathcal{G}$. In other words, $\mathcal{R}$ admits $\mathcal{G}$ as a subset. More precisely, $\mathcal{R}$ is obtained by selecting and then translating edges of $P$ such that no two of them are parallel. In addition, any such selected edge $e$ is shortest in $\langle{e}\rangle_P$. Note that the translation takes place in Lines $8$ and $11$ where $e^\star$ is stored in $\mathcal{R}$ instead of $e$. In this first part of the algorithm, a map $\mu:\mathcal{R}\rightarrow\mathbb{N}$ is also computed in Lines $6$ and $11$ such that $\mu(s)=\left|\langle{s}\rangle_P\right|$ for every segment $s$ in $\mathcal{R}$.

In the second part of algorithm \ref{CZ.computingGZS.algorithm.1}, from Line $15$ to Line $26$, $\mathcal{W}$ is initially set equal to $\mathcal{V}$ and $\mathcal{G}$ to the empty set. The segments in $\mathcal{R}$ that are summands of $P$ are placed in $\mathcal{G}$ in Line $20$ and subtracted from $\mathcal{W}$ by the loop in Lines $21$ to $23$. Lines $18$ and $19$ check whether a segment $s$ in $\mathcal{R}$ is a summand of $P$ using the Minkowski sum of $s$ with the convex hull of $\mathcal{W}$ instead of its Minkowski sum with $P$, allowing for some speedup. This is valid because during the execution of the loop at Line $16$, the $1$-dimensional summands of $P$ remain summands of the convex hull of $\mathcal{W}$ until they are found and subtracted from $\mathcal{W}$.

Let us explain how the subtraction, in Lines $21$ to $23$ of Algorithm \ref{CZ.computingGZS.algorithm.1}, works. By construction, $0$ is a vertex of every segment in $\mathcal{R}$. For any segment $s$ in $\mathcal{R}$, let $s^+$ stand for the non-zero vertex of $s$. If the convex hull of $\mathcal{W}$ admits $s$ as a summand; that is, if it coincides with $R+s$ for some polytope $R$, then $\mathcal{W}$ is naturally partitioned into the points that are vertices of $R$ (because they are the Minkowski sum of a vertex of $R$ with $0$) and the points equal to the sum of a vertex of $R$ with $s^+$. The latter subset is precisely made up of the points in $\mathcal{W}$ that are further displaced by another $s^+$ when the vertex set $\mathcal{X}$ of $\mathrm{conv}(\mathcal{W})+s^+$ is computed in Line~$18$. Therefore, in order to recover the vertex set of $R$, one only needs to subtract $s^+$ from any point in $\mathcal{W}\mathord{\setminus}\mathcal{X}$, which is done in Line $22$, and to keep all the other points in $\mathcal{W}$.

As explained above, all the computations carried out by Algorithm $1$ are polynomial and they are carried out at most a quadratic number of times. In addition, we have seen in Section \ref{CZ.sec.computingZ} that the vertex set of a zonotope can also be computed from its generators in polynomial time.

We therefore obtain the following.

\begin{thm}
The vertex sets of $z(P)$ and $r(P)$ can be computed in polynomial time in $n$, $d$, and the binary size required to store the vertices of $P$.
\end{thm}

\section{Deciding whether a polytope is a zonotope}\label{CZ.sec.decidingZ}

\begin{algorithm}[b]\label{CZ.decidingZ.algorithm.1}

$\mathcal{G}\leftarrow\emptyset$

\For{every subset $\{x,y\}$ of $\mathcal{V}$ such that $x\neq{y}$}{
$e\leftarrow\mathrm{conv}\{x,y\}$

\If{$e$ is an edge of $P$}{
\eIf{some segment $s$ in $\mathcal{G}$ is parallel to $E$}{
\If{$e$ and $s$ have different lengths}{
Return $0$
}
$\mu(s)\leftarrow\mu(s)+1$
}{
$\mathcal{G}\leftarrow\mathcal{G}\cup\{e^\star\}$ and $\mu(e^\star)\leftarrow1$
}
}
}

\For{every segment $s$ in $\mathcal{G}$}{
\If{$\mu_P(s)<{2^{d-1}}$}{Return $0$}
}

\For{every segment $s$ in $\mathcal{G}$}{
Compute the vertex set of $P+s$

\If{the number of vertices of $P$ and $P+s$ is different}{
Return $0$}
}
Return $\mathcal{G}$
\caption{Deciding whether $P$ is a zonotope}
\end{algorithm}

Throughout this section, $P$ is a fixed $d$-dimensional polytope with $n$ vertices, just as in Section \ref{CZ.sec.computingGZS}. Observe that, when $P$ is a zonotope, $z(P)$ is a translate of $P$ and $r(P)$ shrinks to a single point. In particular, Algorithm \ref{CZ.computingGZS.algorithm.1} allows to decide whether a polytope $P$ is a zonotope: this will be the case when the set $\mathcal{W}$ of the vertices of $r(P)$ returned by this algorithm is made up of a single point. In order to solve this decision problem, we can give an alternative algorithm that terminates faster in case the polytope is not a zonotope.

Just as Algorithm~\ref{CZ.computingGZS.algorithm.1}, Algorithm \ref{CZ.decidingZ.algorithm.1} takes as its only input the vertex set $\mathcal{V}$ of $P$. The algorithm returns $0$ when $P$ is not a zonotope. When $P$ is a zonotope, it coincides, up to translation, with $z(P)$. In this case, Algorithm \ref{CZ.decidingZ.algorithm.1} returns the set $\mathcal{G}$ of the generators of $z(P)$. This algorithm consists in three parts. The first part, from Line $1$ to Line $14$ computes a set of line segments that are candidate generators of $P$. Note that this set of line segments is already denoted $\mathcal{G}$ since, if $P$ turns out to be a zonotope, then this set is precisely the set of the generators of  $z(P)$. The computation of $\mathcal{G}$ in Algorithm \ref{CZ.decidingZ.algorithm.1} is very similar to the computation of $\mathcal{R}$ in Algorithm \ref{CZ.computingGZS.algorithm.1}, except that the algorithm immediately terminates in Line $7$ if it finds two parallel edges of $P$ of different lengths. The map $\mu:\mathcal{G}\rightarrow\mathbb{N}$ such that, for any segment $s$ in $\mathcal{G}$, $\mu(s)=\left|\langle{s}\rangle_P\right|$ is computed in Lines $9$ and $11$ like the map $\mu:\mathcal{R}\rightarrow\mathbb{N}$ is in Algorithm \ref{CZ.computingGZS.algorithm.1}.

In the second part of Algorithm \ref{CZ.decidingZ.algorithm.1}, from Line $15$ to Line $19$, every line segment in $\mathcal{G}$ is checked, and the algorithm immediately terminates in Line $17$ if, for such a segment $s$, $\mu(s)<2^{d-1}$. Indeed according to Lemma \ref{CZ.properties.lem.3}, $P$ cannot be a zonotope in this case. In the third part of the algorithm, from Line $20$ to Line $25$, the vertex sets of the Minkowski sums of $P$ with the line segments in $\mathcal{G}$ are computed, and the algorithm terminates in Line $23$ if for such a segment $s$, $P$ and $P+s$ do not have the same number of vertices.

Observe that, if Algorithm \ref{CZ.decidingZ.algorithm.1} does not return $0$ then, for every edge $e$ of $P$, all the segments in $\langle{e}\rangle_P$ are translates of $e$. Moreover, in this case every edge of $P$ is a summand of $P$ because any translate of a summand of $P$ remains a summand of $P$. Therefore, the set $\mathcal{G}$ returned by Algorithm \ref{CZ.decidingZ.algorithm.1} is indeed the set of the generators of $z(P)$. As a consequence, $P$ is a zonotope.

Finally, observe that a zonotope is always centrally-symmetric and, therefore, has an even number of vertices. This very simple test can be done at the beginning of the algorithm to allow for some further speedup.

\bibliography{ComputingZonotopes}
\bibliographystyle{ijmart}

\end{document}